\documentclass[12pt]{amsart}
\usepackage{amsmath,amssymb, amsfonts,amscd,
psfrag,graphicx}

\theoremstyle{plain}
\newtheorem{prop}[equation]{Proposition}

\newtheorem{thm}[equation]{Theorem}
\newtheorem{cor}[equation]{Corollary}
\newtheorem{lem}[equation]{Lemma}

\theoremstyle{definition}

\newtheorem{defn}[equation]{Definition}
\newtheorem{rem}[equation]{Remark}
\newtheorem{exmp}[equation]{Example}
\newtheorem{ques}[equation]{Question}

\numberwithin{equation}{section}

\newcommand{\ZZ}{\mathbb{Z}}
\newcommand{\N}{\mathbb{N}}
\newcommand{\M}{\mathcal{M}}
\newcommand{\F}{\mathcal{F}}
\newcommand{\Lattice}{\mathcal{L}}
\newcommand{\Neigh}{\mathcal{N}}
\newcommand{\LC}{\textrm{LC}}
\newcommand{\NE}{\textrm{NE}}
\newcommand{\Fix}{\textrm{Fix}}
\newcommand{\del}{\textrm{del}}
\newcommand{\lk}{\textrm{lk}}
\newcommand{\lindim}{\textrm{lindim}}
\newcommand{\homdim}{\textrm{homdim}}
\newcommand{\ordim}{\textrm{ordim}}
\newcommand{\lchr}{\textrm{lchr}}
\newcommand{\D}{\Delta}
\newcommand{\G}{\Gamma}

\newcommand{\ka}{\kappa}
\def\maprt#1{\smash{\,\mathop{\longrightarrow}\limits^{#1}\,}}

\begin{document}

\bibliographystyle{plain}

\title{Linear colorings of simplicial complexes and collapsing}

\author{Yusuf Civan and Erg{\" u}n Yal\c c\i n }
\address{Department of Mathematics, Suleyman Demirel
University, Isparta, 32260, Turkey.}

\address{Department of Mathematics, Bilkent
University, Ankara, 06800, Turkey.}

\email{ycivan@fef.sdu.edu.tr \\
yalcine@fen.bilkent.edu.tr}


\thanks{The second author is partially supported by T\" UB\. ITAK-BAYG
through BDP program and by T\" UBA through Young Scientist Award
Program (T\" UBA-GEB\. IP/2005-16)}


\date{\today}

\begin{abstract}
A vertex coloring of a simplicial complex $\Delta$ is called a
\emph{linear coloring} if it satisfies the property that for every
pair of facets $(F_1, F_2)$ of $\Delta$, there exists no pair of
vertices $(v_1, v_2)$ with the same color such that $v_1\in
F_1\backslash F_2$ and $v_2\in F_2\backslash F_1$. We show that
every simplicial complex $\Delta$ which is linearly colored with $k$
colors includes a subcomplex $\Delta '$ with $k$ vertices such that
$\Delta '$ is a strong deformation retract of $\Delta$. We also
prove that this deformation is a nonevasive reduction, in
particular, a collapsing.
\end{abstract}

\maketitle

\section{introduction}
\label{sect:intro}

In this paper, we introduce a notion of linear coloring of a
simplicial complex as a special type of vertex coloring. Recall that
a vertex coloring of an abstract simplicial complex $\D$ with vertex
set $V$ is a surjective map $\ka: V \to [k]$ where $k$ is a positive
integer and $[k]=\{1, \dots , k\}$. We say a vertex coloring is
linear if it satisfies the condition given in the abstract.
Alternatively, a coloring is linear if for every two vertices $u,v$
of $\D$ having the same color, we have either $\F (u) \subseteq \F
(v)$ or $\F (v) \subseteq \F (u)$ where $\F (u)$ and $\F (v)$ denote
the set of facets including $u$ and $v$ respectively. This is
actually equivalent to requiring that the set $\F _i = \{ \F (u) \ |
\ \ka (u)=i \}$ is linearly ordered for every $i\in [k]$, which
explains the rationale for our terminology.

The condition for linear coloring appears naturally when the
multicomplex associated to a colored simplicial complex is studied
closely. For example, in Theorem \ref{thm:isomorphism} we show that
if a simplicial complex is linearly colored then we can recover it
by using the multicomplex associated to it. The multicomplex
associated to a simplicial complex $\D$ is the multicomplex whose
simplices are the color combinations of the simplices on $\D$. We
believe that this association between simplicial complexes and
multicomplexes could be very useful to study the combinatorial
properties of multicomplexes although we do not investigate this
direction in the present work.

Another consequence of requiring a coloring to be a linear coloring
is that it gives us a natural deformation of the colored complex to
a subcomplex of itself where the subcomplex has as many vertices as
the number of colors used. In fact, we can obtain such a deformation
on any subcomplex which satisfies the following condition: Given a
simplicial complex $\D$ and a linear coloring  $\ka$ of $\D$ with
$k$ colors, we call a subcomplex $\D _{\ka} \subseteq \D$ a
\emph{representative subcomplex} if for each $i \in [k]$ there is
one and only one vertex $v$ in $\D _{\ka}$ with $\ka (v)=i$, and if
it has the property that for every pair of vertices $u,v$ with the
same color, we have $\F (u) \subseteq \F (v)$ whenever $u \in \D$
and $v\in \D _{\ka}$. The main result of the paper is the following:

\begin{thm}\label{thm:mainthmintro} Let $\Delta $ be a simplicial
complex on $V$, and let $\kappa : V \to [k]$ be a $k$-linear
coloring map. If $\D _{\ka}$ is a representative subcomplex of $\D$,
then $\D _{\ka}$ is a strong deformation retract of $\D$.
\end{thm}

This allows us to gain information on the homotopy type of a
simplicial complex by coloring it linearly. For example it is clear
that if a simplicial complex can be linearly colored using $k$
colors then its integral (simplical) homology will be zero for
dimensions greater than $k$.

We also introduce the notion of $\LC$-reduction by saying that a
simplical complex $\D$ $\LC$-reduces to its subcomplex $\D '$,
denoted by $\D\searrow_{\LC} \D'$, if there exist a sequence of
subcomplexes $\D =\D _0 \supseteq \D _1 \supseteq \ldots \supseteq
\D_t=\D'$  such that for all $0 \leq r \leq t-1$, the subcomplex $\D
_{r+1}$ is a representative subcomplex of $\D _r $ with respect to
some linear coloring $\ka _r$ of $\D _r$. We study various questions
arising from this definition. For example, we show that if $X_1
\searrow_{\LC} X_2$ and $Y$ is any simplical complex, then $X_1 \ast
Y \searrow_{\LC} X_2 \ast Y$. The main result about $\LC$-reduction
is the following:

\begin{thm}
\label{thm:NEreductionintro} Let $\D$ be a simplicial complex and
$\D'$ be a subcomplex in $\D$. If $\D$ $\LC$-reduces to $\D '$, then
$\D$ $\NE$-reduces to $\D '$ (also called strong collapsing), in
particular $\D$ collapses to $\D'$.
\end{thm}

In fact, Theorem \ref{thm:NEreductionintro} implies Theorem
\ref{thm:mainthmintro}, but we still give a separate proof for
Theorem \ref{thm:mainthmintro} using the basic techniques of poset
homotopy due to Quillen \cite{Qu}. The reason for this is that we
believe that Theorem \ref{thm:mainthmintro} is interesting in its
own right for understanding the topology of simplicial complexes and
should have an independent proof accessible to a topologist. We view
Theorem \ref{thm:NEreductionintro} as a combinatorial version of
Theorem \ref{thm:mainthmintro}.

It turns out that $\LC$-reduction is stronger than the
$\NE$-reduction and hence also stronger than collapsing. In Example
\ref{exmp:irreducible}, we provide an example of a nonevasive
simplicial complex which is not $\LC$-reducible to a point.

In the rest of the paper, we give some applications of
$\LC$-reduction. The first application we give is closely related to
an a theorem by Kozlov \cite{Ko1} about monotone maps and
$\NE$-reduction. We prove that if $\varphi \colon P\to P$ is a
closure operator on a finite poset $P$, then $\D(P)\searrow_{\LC}
\D(\varphi (P) )$, and we conclude that, in this case, $\D(P)$
collapses to $\D(\varphi (P) )$. Our second application is related
to graph coloring. We show that a linear coloring of the
neighborhood complex of a graph gives a (vertex) coloring for the
graph. So, the linear chromatic number of the neighborhood complex
of a simple graph gives an upper bound for the chromatic number of
the graph.

We organize the paper as follows: In Section
\ref{sect:linearcoloring}, we give the definition of a linear
coloring and its equivalent formulations to ease the computations.
Then, in Section \ref{sect:multicomplexes}, we describe an
association between linearly colored simplicial complexes and
multicomplexes. The following three sections contain the  main
results of our work, where we describe the strong deformation of a
simplicial complex induced by a linear coloring, introduce the
notion of $\LC$-reductions, and discuss its connections with known
combinatorial reduction methods such as nonevasive reduction and
collapsing. In particular, we prove Theorem \ref{thm:mainthmintro}
in Section \ref{sect:deformation}, and Theorem
\ref{thm:NEreductionintro} in Section \ref{sect:collapse}.

The last two sections are devoted to applications of
$\LC$-reduction. In Section \ref{sect:posets}, we consider linear
colorings of order complexes of posets and prove the reduction
theorem for closure operators. Finally, in the last section, we
consider the linear colorings of neighborhood complexes associated
to graphs.

\section{Linear coloring of a simplicial complex}
\label{sect:linearcoloring}

We start with some basic definitions related to multisets.

\begin{defn}
\label{defn:multiset} A \emph{multiset} $M$ on a set $S$ is a
function $M\colon S\rightarrow \N :=\{ 0, 1, 2, \dots \}$, where
$M(s)$ is regarded as the number of repetitions of $s\in S$. We say
that $s\in S$ is an \emph{element} of $M$, and write $s\in M$, if
$M(s)>0$. The \emph{cardinality} (or \emph{size}) of a multiset $M$
is defined by $\Arrowvert M \Arrowvert :=\sum_{s\in S} M(s)$.
\end{defn}

Note that every multiset $M$ on $S$ can be regarded as a monomial on
the set $S$ where the degree of $s\in S$ is equal to $M(s)$. The
elements of $M(s)$ will be the elements of $s$ with nonzero degree,
and the cardinality will be equal to the total degree of the
monomial. The usual division relation on monomials gives rise to the
definition of submultisets and the union and the intersection of
multisets can be defined with the following formulas:
\begin{align*}
&(M_1\cup M_2)(s)=M_1(s)+M_2(s); \\
&(M_1\cap M_2)(s)=\textrm{min}(M_1(s),M_2(s)).
\end{align*}

Now we recall the definition of vertex coloring of a simplicial
complex.

\begin{defn}
\label{defn:coloring} Let $\D$ be a finite (abstract) simplicial
complex on $V$. Let $[k]$ denote the set $\{1, \dots, k\}$. A
surjective map $\ka : V \to [k]$ is called a (vertex) coloring of
$\D$ using $k$ colors.
\end{defn}

Given a coloring $\ka$ of a simplicial complex $\D$, we can
associate a multiset to each of its faces as follows: If $S$ is a
face of $\Delta$, then we define the multiset $S_{\ka}$ on $[k]$ by
setting $S_{\ka} (t)$ equal to the order of the set $\{v \in S
\colon \ka (v)=t\}$ for each $t \in [k]$. We define the linear
coloring in its most technical form as follows:

\begin{defn}
\label{defn:linearcoloring} Let $\D$ be a finite abstract simplicial
complex on $V$ and let $\F$ denote the set of all facets of $\D$. A
surjective map $\ka \colon V\rightarrow [k]$ is called a
$k$-\emph{linear coloring} of $\D$ if and only if $\Arrowvert
F_{\ka} \cap F'_{\ka} \Arrowvert =|F\cap F'|$ for any two facets $F,
F'\in \F .$
\end{defn}

Note that if $\D$ is linearly colored with $\ka$, then for distinct
facets $F, F'$ of $\D$, the multisets $F_{\ka}$ and $F' _{\ka}$ must
be also different. Otherwise, we would have $|F\cap F'|=|F|=|F'|$
which cannot happen since $F$ and $F'$ are distinct. We can rephrase
this by saying that the color combinations (with multiplicities)
used in different facets must be different.

Note that every complex with $n$ vertices can be linearly colored
using $n$ colors by giving different color to each vertex. We call a
linear coloring  \emph{trivial} if it is such a coloring.

\begin{defn}
\label{defn:chromaticnumber} The \emph{linear chromatic number} of a
simplicial complex $\D$, denoted by $\lchr (\D)$, is defined to be
the minimum integer $k$ such that $\D$ has a $k$-linear coloring.
\end{defn}

Since there is always the trivial linear coloring, the linear
chromatic number of simplicial complex is well defined and it is
less than or equal to the number of vertices of the complex.

\begin{defn}
\label{defn:type} Let $\D$ be a simplicial complex and let $\ka$ be
a $k$-linear coloring map. Define $V_i:=\{v\in V\;|\;\ka(v)=i\}$ and
set $c^{\ka}_i:=\textrm{card}(V_i)$ for each $i\in [k]$. Then, $\D$
is said to be a linear coloring of \emph{type}
$c_{\ka}(\D)=(c^{\ka}_1,\ldots,c^{\ka}_k)$.
\end{defn}

\begin{exmp}
In Figure \ref{col1}(a), we illustrate a $2$-dimensional simplicial
complex admitting a $2$-linear coloring of type $(3,1)$, whereas
Figure \ref{col1}(b) shows linear coloring of type $(1,1,1,1)$. Note
that the complex in Figure \ref{col1}(b) is a $1$-dimensional
complex with $lchr(\Delta)=4$. For the simplicial complex depicted
in Figures \ref{col1}(c) and \ref{col1}(d), the map given at Figure
\ref{col1}(c) is a $4$-linear coloring of type $(2, 1, 1, 2)$, while
the coloring given in Figure \ref{col1}(d) is not a linear coloring.
\begin{figure}[ht]
\begin{center}
\psfrag{A}{(a)} \psfrag{B}{(b)} \psfrag{C}{(c)} \psfrag{D}{(d)}
\includegraphics[width=4.0in,height=2.2in]{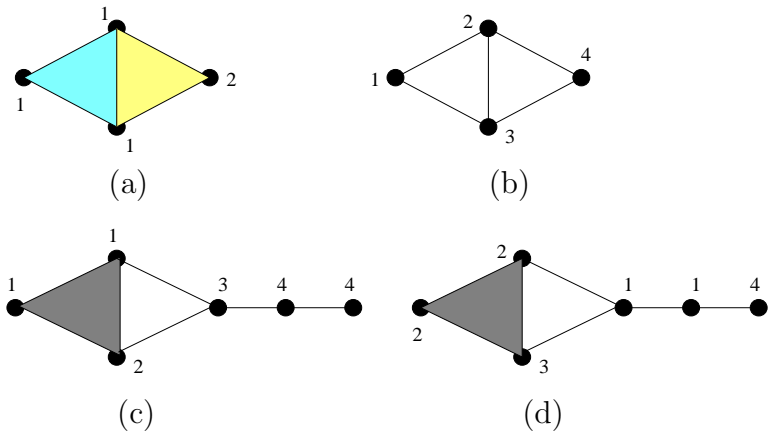}
\end{center}
\caption{Linear colorable complexes and a non-linear coloring}\label{col1}
\end{figure}
\end{exmp}

To understand the definition of linear coloring better, we now give
an equivalent condition for linear coloring. This is the same as the
condition given in the abstract of the paper.

\begin{prop}
\label{prop:eqcond} Let $\Delta$ be a finite abstract simplicial
complex on $V$. A coloring $\kappa\colon V\rightarrow [k]$ of its
vertices is a $k$-\emph{linear coloring} of $\Delta$ if and only if
for every pair of facets $(F_1, F_2)$ of $\Delta$, there exists no
pair of vertices $(v_1, v_2)$ with the same color such that $v_1\in
F_1\backslash F_2$ and $v_2\in F_2\backslash F_1$.
\end{prop}

\begin{proof}
In general $\Arrowvert (F_1)_{\kappa} \cap (F_2)_{\kappa} \Arrowvert
\geq |F_1 \cap F_2 |$ for every pair of facets $(F_1, F_2)$ of
$\Delta$. So, the equality does not hold if and only if there is a
pair of vertices $(v_1, v_2)$ with the same color such that $v_1\in
F_1\backslash F_2$ and $v_2\in F_2\backslash F_1$.
\end{proof}

Note that the above condition for linear coloring can be rephrased
as follows:

\begin{prop}
\label{prop:facetsets}  Let $\Delta$ be a simplicial complex with
vertex set $V$, and let $\kappa\colon V\rightarrow [k]$  be a
coloring of $\D$. For every $v\in V$, let $\F (v)$ denote the set of
facets of $\D$ containing $v$. The coloring $\ka$ is linear  if and
only if for every $i\in [k]$, the set $\F _i= \{ \F (v) : \ka (v)=i
\}$ is linearly ordered by inclusion.
\end{prop}

\begin{proof}
Assume that $\ka $ is a linear coloring. Let $v_1, v_2 \in V$ such
that $\ka (v_1)=\ka (v_2)$. Suppose that there exist facets $F_1 \in
\F (v_1) \backslash \F (v_2)$ and $F_2 \in \F (v_2 )\backslash \F
(v_1)$. Then, it is clear that $v_1 \in F_1 \backslash F_2$ and $v_2
\in F_2 \backslash F_1 $. This contradicts with the fact that $\ka$
is a linear coloring. So, either $\F (v_1) \subseteq \F (v_2)$ or
$\F (v_2 )\subseteq \F (v_1)$ holds. This shows that for each $i$,
the set $\F_i$ is linearly ordered by inclusion. It is clear that
the converse also holds.
\end{proof}

We also have the following observation which will be used later in
the paper.

\begin{prop}\label{prop:extremes}
Let $\Delta$ be a simplicial complex with vertex set $V$, and let
$\kappa\colon V\rightarrow [k]$  be a linear coloring of $\D$. Then,
for each $i\in [k]$, there exists a facet $F$ which includes all the
vertices $v\in V$ with $\ka (v)=i$. On the other extreme, for each
$i \in [k]$, there exists a vertex $v \in V$ such that $v$ lies on
all the facets which include at least one vertex colored with the
color $i$.
\end{prop}

\begin{proof} Take some $i\in [k]$. By Proposition \ref{prop:facetsets},
the set $\F_i=\{ \F (v) : \ka (v)=i \}$ is linearly ordered by
inclusion so there exists a vertex $v \in V$ such that $\ka (v)=i$
and $\F(v)\subseteq \F (u)$ for every $u\in V$ with $\ka (u)=i$. If
we take $F \in \F (v)$, then it is clear that $F$ will include all
the vertices $u \in V$ with $\ka (u)=i$. Note that on the other
extreme, there is a vertex $v\in V$ such that $\ka (v)=i$ and that
$\F(u) \subseteq \F (v)$ for all vertices $u\in V$ with $\ka (u)=i$.
Then, $v$ is included in all the facets which include at least one
vertex colored with the color $i$.
\end{proof}

%%%%%%%%%%%%%%%%%%%%%%%%%%%%%%%%%%%%%%%%%%%%%%%%%%%%%%%%%%%%%%%%%%%%%%%

\section{Multicomplexes associated to linear colorings}
\label{sect:multicomplexes}

In this section, we will discuss an association between
multicomplexes and linearly colored simplicial complexes. We start
with the definition of a multicomplex. More details on this material
can be found in \cite{BV} and \cite{S1}.

\begin{defn}
\label{defn:multicomplex} A \emph{multicomplex} $\Gamma$ is a
collection of multisets over a set $S$ such that if $M\in \Gamma$
and $M'\subseteq M$, then $M'\in \Gamma$. The elements of $\Gamma$
are usually called the {\em faces} of $\Gamma$.
\end{defn}

Note that the faces of $\Gamma$ are ordered by inclusion, giving a
lattice after adjoining a maximal element. We call the resulting
lattice the \emph{face lattice} of $\Gamma$ and denote it by
$L(\Gamma)$. Every multiset $M$ includes a submultiset which is
formed by all its elements with no repetitions. We denote this
submultiset by $u(M)$ and call it the {\em underlying set} of $M$.
If $M$ is a face of a multicomplex $\Gamma$, the underlying set
$u(M)$ of $M$ is called the \emph{underlying face} of $\Gamma$ with
respect to $M$. We have the following simple observation:

\begin{lem}
\label{lem:underlyingsimplex} The collection of all underlying faces
of a multicomplex $\G$ is a simplicial complex. This simplicial
complex is called the {\em underlying simplicial} complex of
$\Gamma$ and denoted by $u(\G)$.
\end{lem}

\begin{proof} Let $S=u(M)$  for some face $M$  of $\Gamma$ and
$S'\subseteq S$. Then $S' \subseteq M$ as a multiset, so $S'$ must
be a face of $\Gamma$. Since $S'=u(S')$, we have $S' \in u(\G)$.
\end{proof}

Now, we consider complexes with a linear coloring.

\begin{prop}
\label{prop:associatedmulticomplex} If $\D$ is a $k$-linear colored
complex with coloring map $\ka$, then the collection $\{ S_{\ka} : S
\in \D \}$ of multisets is a multicomplex.
\end{prop}

\begin{proof} Let $M'$ be a submultiset of a $S_{\ka}$ where $S$
is a simplex in $\D$. Then, it is clear that $S$ has a subset $S'$
such that $S'_{\ka}$ is equal to $M'$.
\end{proof}

\begin{defn}\label{defn:assocmulticomplex} Let $\D$ be a $k$-linear
colored complex with coloring map $\ka$. We call the multicomplex
$\{ S_{\ka} : S \in \D\}$ the {\em associated multicomplex} of the
couple $(\D , \ka)$ and denote it by $\G (\D , \ka)$.
\end{defn}

This gives us an assignment $(\D , \ka ) \to \G (\D , \ka) $ from
the set of linearly colored simplicial complexes to multicomplexes.
The following shows that this assignment is surjective.

\begin{prop}
\label{prop:surjectivity} Given a multicomplex $\G$ over $[k]$,
there exists a simplicial complex $\D$ and a $k$-linear coloring map
$\ka :\D \to [k]$  such that $\G=\G (\D,\ka)$.
\end{prop}

\begin{proof}
Let $\Gamma$ be an arbitrary  multicomplex over $[k]$. For each $i
\in [k]$, let $n_i:=\textrm{max} \{M(i)\colon M\in \Gamma\}$ and let
$V_i:=\{a_r^i\colon 1\leq r\leq n_i \}$. We next define a simplicial
complex $\D(\G)$ on $V:=\cup_{i=1}^k V_i$ as follows: We first
associate a subset $S_M$ of $V$ to every multiset $M\in \G$ by
taking $a_1^i, a_2^i,\ldots, a_j^i\in S_M$ whenever $M(i)=j$ for any
$i\in [k]$. Now, $\D(\G)$ is the $k$-linear colorable simplicial
complex generated by the subsets $F_M\subseteq V$ for which $M$ is a
facet of $\G$, and the linear coloring map $\ka\colon V\rightarrow
[k]$ of $\D(\G)$ is given by $\ka(a_r^i)=i$ for all $i \in [k]$.
\end{proof}

The construction given above gives us a unique simplicial complex
associated to a multicomplex $\G$. Let us denote this simplicial
complex $\D (\G)$. The following shows that the assignment $\G \to
\D (\G)$ is, in fact, inverse to the assignment $(\D , \ka ) \to \G
(\D , \ka)$.

\begin{thm}\label{thm:isomorphism}
Let $\D$ be a simplicial complex on $V$, and let $\ka : V \to [k]$
be a $k$-linear coloring. Suppose $\G=\G (\D , \ka)$ is the
multicomplex associated to the linear coloring $\ka$ and let $\D (\G
)$ be the simplicial complex as in Proposition
\ref{prop:surjectivity}. Then, $\D (\G)$ is isomorphic to $\D$.
\end{thm}

\begin{proof} One can show this using a delicate labeling technique.
Note that the coloring $\ka: V \to [k]$ gives a partitioning of
$V=\cup _{i=1} ^k V_i$ such that $V_i$ is the set of vertices
colored by $i$. Let $n_i$ denote the number of elements in $V_i$ for
each $i \in [k]$. As before let $\F (v) $ denote the set of facets
in $\D$ including $v$ as a vertex. Recall that by Proposition
\ref{prop:facetsets}, for each $i\in [k] $, the set $\F _i =\{ \F (v
) : v \in V_i\}$ is linearly ordered by inclusion. We can label the
vertices of $\D$ in the following way: Let $V= \{v_{r}^i \colon i
\in [k], r \in [n_i] \}$ where for all $i$, the vertex $v_{r}^i  $
belongs to $V_i$ and  $\F (v_t ^i) \subseteq \F (v_{r}^i)$ whenever
$1 \leq r\leq t \leq n_i$.

Recall that the simplicial complex $\D(\G)$ on $V:=\cup_{i=1}^k V_i$
is defined as follows. The subset $S_M$ of $V$ to every multiset
$M\in \G$ is defined by taking $a_1^i, a_2^i,\ldots, a_j^i\in S_M$
whenever $M(i)=j$ for any $i\in [k]$. Now, $\D(\G)$ is the
simplicial complex generated by the subsets $F_M\subseteq V$ for
which $M$ is a facet of $\G$.

We claim that the assignment $f: \D \to \D (\G)$ defined by $f (v_r
^i )= a_r ^i$ for every $i \in [k]$ and $r\in [n_i]$ is an
isomorphism of simplicial complexes. To prove this claim, it is
enough to show that $S$ is a simplex in $\D$ if and only if $f(S)$
is a simplex in $\D (\G)$. Note that we can prove each direction
starting with a facet. Let $F$ be a facet in $\D$. To show that
$f(F)$ is a simplex in $\D (\G )$, we need to show that $F$
satisfies the property that if $v^i _t \in F$, then $v ^i _r$ is in
$F$ for every $1 \leq r \leq t$. This follows from the fact that $\F
(v_t ^i) \subseteq \F (v_{r}^i)$ for every $1 \leq r\leq t \leq
n_i$. So, $f(F) \in \D (\G)$ as desired. For the other direction,
let $F$ be a facet in $\D (\G)$, and let $M$ be the corresponding
face in $\G$. Then, there is a facet $F'$ in $\D$ such that for each
$i \in [k]$, a vertex from $V_i$ appears exactly $M(i)$ times.
Recall that the facets of $\D$ satisfy the property that if $v^i _t
$ is in a facet, then $v ^i _r$ is also in that facet for every $1
\leq r \leq t$. So, we can conclude that $F'=f^{-1}(F)$, and hence
$f^{-1} (F)$ is in $\D$. This completes the proof.
\end{proof}

This shows, in particular, that we can recover a linearly colored
simplicial complex from its associated multicomplex. Another way to
state the above result is the following:

\begin{cor}\label{cor:labeling} Let $\Delta $ be a simplicial
complex on $V$, and let $\kappa : V \to [k]$ be a $k$-linear
coloring. For each $i \in [k]$, let $V_i:=\{v \in V \colon \ka (v)=i
\}$ and let $n_i=|V_i|$. Then, we can label the vertices of $\D$ in
such a way that $V= \{v_{r}^i \colon i \in [k], r \in [n_i] \}$ and
that whenever $v^i _t $ is in a facet of $\D$, then $v ^i _r$ is
also in that facet for every $1 \leq r \leq t$.
\end{cor}

The labeling technique given in the above corollary can also be used
to produce some poset maps between the face posets of the simplex
$\D$, the associated multicomplex $\Gamma (\D, \ka )$, and the
underlying simplicial complex $u (\Gamma (\D , \ka ))$. We now
explain these.

\begin{lem}\label{lem:coloringposetmap}
Let $\Delta $ be a simplicial complex linearly colored with a
coloring map $\kappa : V \to [k]$, and let $\Gamma$ be the
associated multicomplex. Then the map $$ c: \Delta \to \Gamma $$
defined by $S \to S_\kappa$ for every $S \in \Delta$ is a poset map
(between corresponding face posets).
\end{lem}

\begin{proof} This is clear since for every $S' \subseteq S$, the
number of times a color used in $S'$ is less than or equal to the
number of times it is used in $S$.
\end{proof}

We also have the following:

\begin{lem}\label{lem:underlyingposetmap} Let  $\Gamma$ be a
multicomplex, and let $u(\Gamma)$ denote its underlying simplicial
complex. The canonical map $$u: \Gamma \to u(\Gamma )$$ defined by
$M\to u(M)$ for every $M\in \Gamma $ is a poset map.
\end{lem}

\begin{proof} If $M'\leq M$, then $M'(t)\leq M(t)$ for all $t$. In
particular, if $M'(t)> 0$, then $M(t)>0$.
\end{proof}

Given a linear coloring $\kappa : V( \D ) \to [k]$ for $\D$, let
$$\varphi _{\ka} : \Delta \to u(\Gamma (\D , \kappa ))$$
denote the composite map $u \circ c$. It is clear by the above two
lemmas that $\varphi _{\ka}$ is a poset map between face posets of
two simplicial complexes. So, considered as a map between two
simplicial complexes, it is a simplicial map.

\begin{prop}\label{prop:imbedding}
Let $\D$ be a simplicial complex and $\ka : V \to [k]$ be a
$k$-linear coloring. Then, there exists a simplicial map
$$i_{\ka}: u(\G (\D , \ka)) \to \D$$ such that
$\varphi _{\ka} \circ i_{\ka}=id $.
\end{prop}

\begin{proof} Suppose that the vertices of $\D$ are labeled as in
Corollary \ref{cor:labeling}. We first show that for each $S \in u
(\G (\D, \ka ))$, the set $V_S=\{ v_{1}^s : s \in S \}$ is a simplex
in $\D$. Note that if $F$ is a facet of $\D$ such that $v^s _i \in
F$ for some $s \in S$ and some $i\in [n_i]$, then the vertex $v_1
^s$ belongs to $F$. Since $S$ is equal to $u(M)$ for some multiset
$M$ in $\G =\G (\D, \ka )$, the set $S$ considered as a multiset
belongs to multicomplex $\G$. We also observed earlier that there is
a one-to-one correspondence between facets of $\D$ and $\Gamma$, so
we can choose a facet $F$ of $\D$ such that $S \leq F_{\ka}$. This
facet has to include $v^s _1$ for all $s \in S$ by the above
argument. So, $V_S$ is a simplex of $\D$.

Let $i_{\ka}: u(\G (\D , \ka)) \to \D$ be the map defined by
$i_{\ka} (S) =V_S$ for every $S \in u(\G (\D , \ka ))$.  It is easy
to see that $i_{\ka}$  satisfies the desired properties.
\end{proof}

Note that the simplicial map $i_{\ka}$ is not uniquely defined in
general. This is because the set of faces $\F_i = \{ \F (v ) : \ka
(v ) =i \}$ can be linearly ordered in many different ways, and as a
result of these different orderings there could be more than one
vertex that we can choose as the vertex with label $v^i _1$. On the
other hand the subcomplexes which can be the image of $i_{\ka}$ have
something in common.  Their vertices are colored with distinct
colors and have the property that for every pair of vertices $x,y$
with $x \in \D$ and $y\in i _{\ka} (u(\G (\D , \ka))) $ such that
$\ka (x)=\ka (y)$, we have $\F (x) \subseteq \F (y)$. Conversely any
subcomplex having these properties is the image of $i_{\ka}$ for
some choice of ordering. We will study such subcomplexes further in
the next section.

\section{Deformation to a representative subcomplex}
\label{sect:deformation}

In this section we prove Theorem \ref{thm:mainthmintro} stated in
the introduction. We start with the definition of representative
subcomplex.

\begin{defn}
\label{defn:representative} Let $\D$ be a simplicial complex with
linear coloring map $\ka : V(\D) \to [k]$ where $V(\D)$ denotes the
vertex set of $\D$. A subcomplex $\D _{\ka}$ of $\D$  is said to be
a \emph{representative subcomplex} with respect to $\ka$ if for each
$i \in [k]$ there is one and only one vertex in $x \in \D _{\ka}$
with $\ka (x)=i$ and if it has the property that for every pair of
vertices $x,y$ with $x \in \D$, $y\in \D _{\ka}$ and $\ka (x)=\ka
(y)$, we have $\F (x) \subseteq \F (y)$, where $\F (x)$ and $\F (y)$
denote the set of facets including $x$ and $y$ respectively.
\end{defn}

Although a linearly colored complex may have many different representing
subcomplexes, the following result shows that as simplicial
complexes they are all same.

\begin{prop}
\label{prop:uniquenessofinc} Let $\D$ be a simplicial complex with
linear coloring $\ka$. Suppose that $\D _{\ka}$ and $\D _{\ka } '$
are two subcomplexes of $\D$ which are representative with respect
to $\ka$. Then, $\D _{\ka}$ and $\D _{\ka } '$ are isomorphic as
simplicial complexes.
\end{prop}

\begin{proof} Let $x,y $ be two vertices in a simplicial complex
with $\F (x)=\F (y)$. Consider the map $f: V(\D) \to V(\D)$ such
that $f(x)=y, f(y)=x$ and $f(z)=z$ for all the other vertices. We
claim that $f$ extends to an isomorphism of simplicial complexes.
For this it is enough to show that if $S \in \D$, then $f(S)\in \D$.
This is clear if $x,y$ are both in $S$ or if neither of them are in
$S$. Suppose $S$ is such that $x\in S$ and $y \not \in S$. Let $F$
be a facet that includes $S$. Since $x\in F$, we must have $y \in F$
by the assumption that $\F (x)=\F (y)$. This gives that $f(F)=F$.
From this we can conclude that $f(S) \subseteq F$ and hence $f(S)$
is a simplex in $\D$. Similarly, if $S$ is a simplex with $y \in S$
and $x \not \in S$, we can prove again $f(S)$ is in $\D$ using the
equality $\F (x)=\F (y)$.

Let $\D_{\ka}$ and $\D '_{\ka}$ be two different choices of
representative subcomplexes. Composing isomorphisms of the above
type, we can find an isomorphism $f:\D \to \D$ such that  $f$ takes
the image of $\D _{\ka}$ to the image of $\D '_{\ka}$.
\end{proof}

We are now ready to prove Theorem \ref{thm:mainthmintro}:

\begin{proof}[Proof of Theorem \ref{thm:mainthmintro}]
We need to show that the composition
$$f: \D \maprt{i_{\ka}\circ\varphi _{\ka}} \D _{\ka} \maprt{inc}
\D$$ is homotopic to identity with a homotopy relative to $\Delta
_{\ka}$. Note that there exists a unique inclusion $i_{\ka}$ once
$\D _{\ka}$ is chosen. Also, it is clear that $f$ is a poset map
between corresponding face posets. If there exists another poset map
$g: \D \to \D $ such that $ S \leq g(S)\geq f(S)$ for all $S\in \D$,
then by Quillen's criteria for homotopy equivalence of poset maps
(see, for example, \cite{Qu}), we can conclude that $id\simeq
g\simeq f$. Below we show that for every $S\in \D$, the set $S \cup
f(S)$ is a simplex of $\D$. This allows us to define $g:\D \to \D $
as the map $g(S)= S \cup f(S)$ and conclude that $f$ is homotopic to
identity. Since both $f$ and $g$ are equal to identity on $\Delta
_{\ka}$, the required relativeness condition for the homotopy also
holds.

To show that $S \cup f(S)$ belongs to $\D$ for all $S\in\D$, we use
the labeling given in Corollary \ref{cor:labeling}. Suppose that the
vertices of $\D$ are labeled as in Corollary \ref{cor:labeling}.
Note that $$f(S)= \{ v_1 ^i : i \in u(S_{\ka}) \}$$ for every
simplex $S\in \D$. Let $S$ be a simplex in $\D$ and $F$ be a facet
including $S$. If the color $i$ is used to color a vertex in $S$,
then $S$ must include a vertex of the form $v^i _r$ for some $r\in
[n_i]$. The way we have chosen the labeling implies that $v^i _1 \in
F$. Since this is true for all $i \in u(S_{\ka})$, we can conclude
that $f(S) \subseteq F$. Since $F$ includes both $S$ and $f(S)$, it
includes $S\cup f(S)$. This shows that $S \cup f(S)$ is a simplex of
$\D$. This completes the proof.
\end{proof}

The following is an immediate corollary of  Theorem
\ref{thm:mainthmintro}.

\begin{cor}\label{cor:homology}
Let $\Delta $ be a $k$-linear colorable simplicial complex. Then, $H_i
(\D , \ZZ )=0$ for all $i \geq k$.
\end{cor}

Another important consequence of Theorem \ref{thm:mainthmintro} is
that it provides a lower bound for the linear chromatic number of a
simplicial complex by the topology of the complex. To see this, we
first introduce some terminology about connectedness. Let
$\widetilde{H}_i(\D)$ denote the \emph{reduced simplicial homology
groups} (over $\ZZ$) of a simplicial complex $\D$ (see \cite{Mu} for
details). A simplicial complex $\D$ is said to be $k$-\emph{acyclic}
if $\widetilde{H}_r(\D)=0$ for all $r\leq k$, and it is called
\emph{acyclic} if it is $k$-acyclic for all $k\in \ZZ$. Further,
$\D$ is called $k$-\emph{connected} if it is $k$-acyclic and simply
connected, $k\geq 1$. The following is the linear coloring analogue
of a well-known result of Lov\'{a}sz on graph colorability (see
\cite{Lo}).

\begin{cor}
\label{cor:k-connected} If $\D$ is non-acyclic and $k$-connected
($k\geq 1$), then $\lchr(\D)\geq k+3$.
\end{cor}

\begin{proof}
Assume that $\D$ admits a $(k+2)$-linear coloring $\ka$ and let $\D
_{\ka}$ be a representative subcomplex of $\D$ with respect to
$\ka$. Then, $\D$ is homotopy equivalent to $\D_{\ka}$ by Theorem
\ref{thm:mainthmintro}, where $\D_{\ka}$ is a simplicial complex
with $k+2$ vertices. Such a complex is at most $(k+1)$-dimensional.
Since $\D$ is non-acyclic, the dimension of $\D_{\ka}$ cannot be
less than $k+1$ by $k$-connectivity. On the other hand, if
$\textrm{dim}(\D_{\ka})=k+1$, then it is a $(k+1)$-simplex which is
contractible; hence, it is acyclic, a contradiction.
\end{proof}

%%%%%%%%%%%%%%%%%%%%%%%%%%%%%%%%%%%%%%%%%%%%%%%%%%%%%%%%%%%%%%%%%%%%%%%%%%%%%%%%%%%

\section{$\LC$-reduction of a simplicial complex}
\label{sect:LC-reduction}

In this section we introduce the concept of $\LC$-reduction and
study its basic properties. We start with the definition of
$\LC$-reduction.

\begin{defn}\label{defn:LC-reduction}
Let $\D$ be a simplicial complex and $\D'$ be a subcomplex of
$\Delta$. If there exist a sequence of subcomplexes $\D =\D _0
\supseteq \D _1 \supseteq \ldots \supseteq \D_t=\D'$ such that $\D
_{r+1}$ is a representative subcomplex in $\D _r $ with respect to
some linear coloring $\ka _r$ of $\D _r$ for all $0\leq r \leq t-1$,
then we say $\D$ $\LC$-\emph{reduces} to $\D'$, and write
$\D\searrow_{\LC} \D'$.
\end{defn}

By Theorem \ref{thm:mainthmintro}, it is easy to see that if $\D$
$\LC$-reduces to a subcomplex $\D'$, then $\D'$ is a strong
deformation retract of $\D$.

For our purposes it is desirable to be able to express an
$\LC$-reduction as a composition of $\LC$-reductions which are
primitive in some sense. In this context, the appropriate definition
of  primitiveness can be given as follows:

\begin{defn}\label{defn:primitive}
A linear coloring of a simplicial complex $\D$ with $n$ vertices is
called a \emph{primitive linear coloring} if it uses exactly $n-1$
colors. An $\LC$-reduction is called primitive if it involves only
one linear coloring and that coloring is primitive.
\end{defn}

Note that if $\ka$ is a primitive linear coloring then there is a
pair of vertices $u,v$ in $\D$ such that $\ka (u)=\ka (v)$ and the
remaining vertices of $\D$ are colored using distinct colors. By the
condition of a linear coloring, we have either $\F (u) \subseteq \F
(v)$ or $\F (v) \subseteq \F (u)$. In the first case, the subcomplex
$\del _{\D} (u)= \{ S \in \D \ | \ u \not \in S \}$ will be a
representative subcomplex, and in the second case $\del _{\D} (v)=
\{ S \in \D \ | \ v \not \in S \}$ will be representative. In the
case of equality either of these sets can be taken as a
representative subcomplex. Note that an $\LC$-reduction $\D \searrow
_{\LC} \D'$ is primitive if and only if the number of vertices in
$\D'$ is exactly one less than the number of vertices in $\D$.

\begin{prop}\label{prop:primitive}
Any $\LC$-reduction $\D\searrow_{\LC} \D '$ can
be expressed as a sequence of primitive $\LC$-reductions.
\end{prop}

\begin{proof} It is enough to prove the proposition for a
$\LC$-reduction involving only one coloring. So, we can assume  $\D
'= \D _{\ka}$ for some coloring $\ka$ of $\D$. Suppose that the
vertices $\D$ are labeled as in Corollary \ref{cor:labeling}. So, if
$V$ is the set of vertices of $\D$, then we can write $V= \{ v_{r}^i
\colon i \in [k], r \in [n_i] \}$ where  $\F (v_t ^i) \subseteq \F
(v_r ^i)$ whenever $1 \leq r\leq t \leq n_i$. We can assume that $\D
_{\ka }$ is the subcomplex generated by the vertices $\{ v^i _1 \ |
\ i =1, \dots , k \}$.

Let $\ka (i,j)$ denote the primitive linear coloring involving
vertices $v^i _j$ and $v^i _{j+1}$ for $i=1, \dots k$ and $j=1,
\dots , n_i -1$. It is easy to see that if we apply $\LC$-reductions
associated to primitive linear colorings $ \ka (i, n_i-1), \ka (i,
n_i-2), \dots , \ka (i, 1)$ in this order for each $i=1,\dots k$,
then we obtain an $\LC$-reduction to $\D _{\ka}$.
\end{proof}

Some complexes cannot be $\LC$-reduced further to any proper
subcomplex.

\begin{defn}\label{defn:irreducible}
A simplicial complex $\D$ on a set $V$ is called
$\LC$-\emph{irreducible} if it admits only a trivial linear
coloring.
\end{defn}

The following is clear from the definition.

\begin{prop}
\label{prop:irreducibilitycondition}
A simplicial complex $\D$ is $\LC$-irreducible if and only if for every pair of
vertices $u, v$, the facet sets $\F (u)$ and $\F(v)$ are not comparable
by inclusion.
\end{prop}

A typical example of an $\LC$-irreducible complex is the boundary of a simplex.
Another example would be a complex whose realization is an $n$-gon.

It is easy to see that every simplicial complex $\D$ $\LC$-reduces
to an $\LC$-irreducible subcomplex, although  the resulting
$\LC$-irreducible subcomplex can be quite different depending on the
choices we make. Let us call a subcomplex $\D'$ of $\D$ an
$\LC$-\emph{core} of $\D$ if it is irreducible and if $\D$
$\LC$-reduces to it. The homotopy type of an $\LC$-core is uniquely
determined by the homotopy type of $\D$, but it is not easy to see
what other properties of $\LC$-cores of $\D$ are invariants of $\D$.
It is reasonable to ask:

\begin{ques}\label{ques:uniquenessofcore}
Let $\D$ be a simplicial complex and $\D _1$ and $\D _2$
are two different $\LC$-cores for $\D$. Is it true that $\D _1 $
and $\D _2$ are isomorphic as simplicial complexes?
\end{ques}

At this point we do not know the answer to this question. One would
expect that at least the number of vertices of a core is an
invariant of the simplicial complex. Until finding an answer to this
question we can define such an invariant as follows:

\begin{defn}\label{defn:lindim} Let $\D$ be a simplicial complex.
The \emph{linear dimension} of $\D$, denoted by $\lindim (\D)$, is
defined to be the smallest integer $n$ such that $\D$ has a core
with $n$ vertices.
\end{defn}

Note that $\lindim (\D)$ is also the smallest integer $n$ such that
$\D$ $\LC$-reduces to a simplicial complex with $n$ vertices. It is
easy to see that linear dimension is related to the homological
dimension of the complex. Recall that the \emph{homology dimension}
$\homdim (\D)$ of a finite simplicial complex $\D$ is defined to be
the integer
$$\homdim (\D):=\mathrm{min}\{ i\ | \ \widetilde{H}_j(\D; \ZZ )=0\;
\mathrm{for\;all}\;j>i \}$$ with the convention that
$\widetilde{H}_{-1}(\D; \ZZ )=\ZZ$. We can easily adopt the proof of
Corollary \ref{cor:k-connected} to obtain the following.

\begin{prop}\label{prop:lhomdim}
For any finite simplicial complex
$\D$, we have $$\lchr (\D) \geq \lindim(\D)\geq \homdim(\D)+2.$$
\end{prop}

An interesting family of simplicial complexes are the ones with
linear dimension equal to one. These are the complexes which can be
$\LC$-reduced to a point. We say a simplicial complex $\D$ is
\emph{$\LC$-contractible} if $\D \searrow _{\LC} \{x\}$ for some
vertex $x$ of $\D$. We use this terminology later in the paper.

Now, we investigate the behavior of $\LC$-reduction under the join
operator. Recall that the join of two simplicial complexes $X$ and
$Y$, denoted by $X \ast Y$, is defined as the simplicial complex
which includes both $X$ and $Y$ as subcomplexes and includes also
the sets of the form $S \cup T$ where $S \in X$ and $T\in Y$.

\begin{prop}
\label{prop:join}
Let $X_1 \searrow _{\LC} X_2$ and let $Y$ be an arbitrary simplicial complex.
Then,  $X_1 \ast Y \searrow _{\LC} X_2 \ast Y $.
\end{prop}

\begin{proof} It is enough to prove the result for a primitive $\LC$-reduction.
Let $X_1 \searrow _{\LC} X_2$ be a primitive reduction involving
vertices $u,v \in X_1$. Without loss of generality we can assume $v
\in X_2$. Recall that in this case $X_2$ is the subcomplex $\del
_{X_1} (u)=\{ S\in X_1 \ | \ u \not \in S \}$. Since $\del _{X_1
\ast Y} (u)= \del _{X_1} (u)\ast Y$, we just need to show that
primitive coloring involving $u$ and $v$ is still a linear coloring
in $X_1 \ast Y$. We know that $\F (u) \subseteq \F (v)$ in $X_1$.
Let $F$ be a facet of $X_1 \ast Y$ including the vertex $u$. Then
either $F$ is a facet of $X$ or $F$ is of the form $S \cup T$ where
$S$ and $T$ are facets of $X$ and $Y$ respectively. In the first
case, $F \in \F(u)$, so $v \in F$ can be seen easily. In the second
case, the facet $S $ belongs to the set $\F (u)$, and again we can
conclude $v \in S$. This gives $v \in F$ since $F=S \cup T$. This
shows that the inclusion $\F (u) \subseteq \F (v)$ still holds for
facet sets in $X_1 \ast Y$. This completes the proof.
\end{proof}

%%%%%%%%%%%%%%%%%%%%%%%%%%%%%%%%%%%%%%%%%%%%%%%%%%%%%%%%%%%%%%%%%%%%%%%

\section{$\LC$-reduction, nonevasive reduction and collapsing}
\label{sect:collapse}

The aim of this section is to prove Theorem
\ref{thm:NEreductionintro} stated in the introduction. We first
recall the definition of collapsing.

\begin{defn}\label{defn:collapse}
A face $S$ of a simplicial complex $\D$ is called \emph{free} if $S$
is not maximal and there is a unique maximal face in $\D$ that
contains $S$. If $S$ is a free face of $\D$ then the simplicial
complex $\D[S]:=\D\backslash \{T \in \D \ | \ S \subseteq T \}$ is
called an \emph{elementary collapse} of $\D$. If $\D$ can be reduced
to a subcomplex $\D'$  by a sequence of elementary collapses, then
we say $\D$ \emph{collapses} to $\D'$. In this case, we write
$\D\searrow \D'$. If a complex collapses to a point then we say it
is \emph{collapsible}.
\end{defn}

We start with the following result:

\begin{prop}
\label{prop:collapsibility} Let $\D$ be a simplicial complex and
$\D'$ be a subcomplex in $\D$. If $\D \searrow _{\LC} \D '$, then
$\D \searrow \D'$.
\end{prop}

\begin{proof} It is enough to prove the proposition for a primitive
linear coloring. So, assume that $\D ' =\D _\ka $ for some primitive
linear coloring $\ka$ which involves vertices $u$ and $v$.  Without
loss of generality, we may assume that $u$ lies on $\D _{\ka}$. Note
that this implies in particular that $\F(v)\subseteq \F(u)$. Let
$F_1\in \F(v)$ be given. Then, we claim that the face
$S_1:=F_1\backslash \{u\}$ is contained only in $F_1$, i.e., it is
free in $\D$. Indeed, if $F'$ is any facet containing $S_1$, then
$v\in F'$. This gives $u\in F'$ because $\F(v)\subseteq \F(u)$. But
then $F_1 \subseteq F'$, and we can conclude that $F_1=F'$. Let
$\D_1$ denote the elementary collapse of $\D$ through the face
$S_1$, that is, $\D_1=\D[S_1]$. For the simplicial complex $\D_1$,
we note that any facet containing the vertex $v$ must also contain
$u$. Therefore, we may similarly collapse $\D_1$ by choosing a facet
$F_2 $ of $\D_1$ containing $v$. We iterate the same process until
we obtain a simplicial complex $\D_m$ in which $\F (v)$ is empty. It
is easy to see that $\D_{m}=\del _{\D} (v)$, and hence it is equal
to $\D_{\ka}$.
\end{proof}

The converse of Proposition \ref{prop:collapsibility} does not hold
in general.
\begin{exmp}\label{exmp:irreducible}
Let $\D$ be the $2$-dimensional simplicial complex on $V=\{a,b,c,d,e,f\}$
with the set of facets
$$\F(\D)=\{\{a,b,c\},\{a,b,e\},\{a,d,e\}, \{b,e,f\},
\{d,e,f\},\{b,c,f\}, \{c,d,f\}\}.$$ The realization of $\D$ is given
in Figure \ref{col5}, where the picture is intended to be three
dimensional like a pyramid. Note also that the interior of the
shaded simplex is not part of the complex. It is clear that $\D$ is
collapsible and $\NE$-reduces to a point (i.e. nonevasive), but it
does not $\LC$-reduce to a point (in fact it is $LC$-irreducible).
\begin{figure}[ht]
\begin{center}
\includegraphics[width=2.0in,height=1.4in]{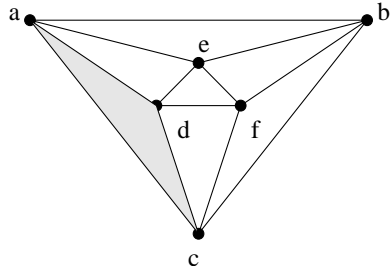}
\end{center}
\caption{A collapsible $\LC$-irreducible simplicial
complex}\label{col5}
\end{figure}
\end{exmp}

Another type of reduction of simplicial complexes is nonevasive
reduction (see Kozlov \cite{Ko1}, Welker \cite{We}) which is also
known as strong collapsing (see Kurzweil \cite{Ku}). Recall that for
a vertex $v$ in a simplicial complex $\D$, the deletion of $v$ is
defined as the subcomplex $\del_{\D} (v)=\{ S \in \D \ | \ v\not \in
S \}$ and the link of $v$ in $\D$ is defined as the subcomplex $\lk
_{\D} (v) = \{ S  \in \D \ | v \not \in S , S\cup \{ v \} \in \D \}
$. Nonevasiveness of a simplicial complex is defined inductively by
declaring that a point is nonevasive and a simplicial complex $\D$
is nonevasive if it has a vertex $v$ such that both its deletion
$\del _{\D} (v)$ and its link $\lk _{\D} (v)$ are nonevasive.

\begin{defn}\label{defn:nonevasivereduction}
Let $\D$ be a simplicial complex and $\D'$ be a subcomplex of $\D$.
We say that $\D$ $\NE$-\emph{reduces} to $\D'$, denoted by
$\D\searrow_{\NE} \D'$, if there exist a sequence $\D =\D^1,\D^2,
\ldots, \D^{t+1}=\D'$ of subcomplexes and a sequence of vertices
$v_1, \dots v_t$ such that $V(\D^{r})=V(\D^{r+1}) \cup \{v_r\} $ and
$\lk _{\D _r} (v_r ) $ is nonevasive for any $1\leq r\leq t$.
\end{defn}

We have the following:

\begin{prop}
\label{prop:NEreduction} Let $\D$ be a simplicial complex and $\D'$
be a subcomplex of $\D$. If $\D \searrow _{\LC} \D '$, then $\D
\searrow _{\NE} \D'$.
\end{prop}

\begin{proof} As before it is enough to prove the proposition for a
primitive linear coloring. Let $\D'=\D_{\ka}$ where $\ka$ is a
primitive coloring involving vertices $u$ and $v$. Without loss of
generality we can assume $u \in \D '$. We have $\F(v) \subseteq
\F(u)$. We claim that $\lk _{\D} (v)$ is nonevasive. This will imply
that $\D \searrow _{\NE} \D'$  as desired.

It is well known that if a simplicial complex is a cone then it is
nonevasive. So, it is enough to show that $\lk _{\D} (v)$ is a cone.
Let $S$ be a simplex in $\lk _{\D} (v)$. Then, $S\cup \{v \}$ is a
simplex in $\D$. Let $F$ be a facet of $\D$ which includes $S\cup
\{v \}$. Since $F \in \F (v)$, we have $F \in \F (u)$ by our
assumption that $\F(v) \subseteq \F(u)$. This implies that $S\cup \{
u\}$ is a simplex in $\lk _{\D} (v)$. We have shown that for every
simplex $S$ in $\lk _{\D} (v)$, $S\cup \{ u\}$ is also a simplex in
$\lk _{\D} (v)$. This means $\lk _{\D} (v)$ is a cone with cone
point $u$.
\end{proof}

Note that the proof of Theorem \ref{thm:NEreductionintro} is now
complete.

\begin{rem}  It is well known that nonevasive reduction
is a collapsing by a result of Kahn, Saks, and Sturtevant (see
Proposition 1 in \cite{KSS}). So, Proposition
\ref{prop:collapsibility} can also be obtained as a corollary of
Proposition \ref{prop:NEreduction}.
\end{rem}

%%%%%%%%%%%%%%%%%%%%%%%%%%%%%%%%%%%%%%%%%%%%%%%%%%%%%%%%%%%%%%%%%%%%%%%%%%%%%%%

\section{Linear coloring of posets}
\label{sect:posets}

Let $P$ be a finite partially ordered set. We denote by $\D(P)$ its
order complex, i.e., the set of all chains in $P$. When $P$ has
maximal and minimal elements, we denote them by $\hat{0}$ and
$\hat{1}$, respectively. The elements of $P$ that cover $\hat{0}$
are called \emph{atoms}, and the elements that are covered by
$\hat{1}$ are called \emph{coatoms}. We denote the set of atoms and
coatoms of a bounded poset $P$ by $at(P)$ and $co(P)$ respectively.
We write $\overline{P}$ for the poset $P\backslash \{\hat{0},
\hat{1}\}$, and call it the proper part of $P$. The set of maximal
chains of $P$ is denote by $\M$, and in particular $\M_x$ denotes
the maximal chains containing the element $x\in P$. For a given
subset $S\subseteq P$, we denote by $\bigwedge S$ and $\bigvee S$,
the greatest lower bound and the least upper bound (when exist) of
$S$ respectively.

Throughout, by a linear coloring of $P$, we  mean a linear coloring
of $\D(P)$. We may rephrase the definition of a linear coloring for
posets as follows.

\begin{lem}\label{lem:multiposet}
A surjective mapping $\ka\colon P\to [k]$ is a $k$-linear coloring
of $P$ if and only if $\ka(x)=\ka(y)$ implies either $\M_x \subseteq
\M_y$ or $\M_y \subseteq \M_x$ for any two elements $x,y \in P$.
\end{lem}

This implies, in particular, that in a linearly colored poset $P$
any two elements $x,y\in P$ having the same color must be
comparable. In fact, more is true. Let $P$ be a poset linearly
colored with $\ka$, and let $x,y \in P$ be such that
$\ka(x)=\ka(y)$. Suppose $\M_x \subseteq \M_y$. Let $z$ be an
element in $P$ such that $x$ is comparable with $z$, i.e, either
$x<z$ or $z< x$. Then, there is a maximal chain $M$ including $x$
and $z$. Since $\M _x$ is included in $\M _y $, the chain $M$ must
also include $y$. Thus, $z$ and $y$ are also comparable. Similarly,
we can show that if $\M_x\subseteq \M_y$, then every element of $P$
which is comparable with $y$ is also comparable with $x$. We define
the following:

\begin{defn}\label{defn:domination}
Let $P$ be a poset and $x, y \in P$. We say $y$ \emph{dominates}
$x$, denoted by $x \prec  y$,  if every element $z$ which is
comparable with $x$ is also comparable with $y$.
\end{defn}

We have seen above that in a linearly colored poset $P$ any two
elements $x,y\in P$ having the same color must be comparable by
domination. The converse of this statement also holds:

\begin{prop}\label{prop:domination}
Let $P$ be a poset and $\ka : P \to [k]$ be a coloring of $P$. Then,
$\ka$ is a linear coloring if and only if for every pair $x,y \in P$
with $\ka (x)=\ka (y)$, either $x \prec y$ or $y \prec x$.
\end{prop}

\begin{proof} We only need to prove one direction. Let $x,y \in P$
be such that $\ka (x)=\ka (y)$ and $x\prec y$. Then every element $z
\in P$ which is comparable with $x$ is also comparable with $y$. We
claim that in this case the inclusion $\M_x\subseteq \M_y$ holds.
Let $M$ be a maximal chain in $\M _x$. Note that all the elements in
$M$ are comparable with $x$, so they must be also comparable with
$y$. If $y$ is not in $M$, then by adding $y$ to $M$ we would get a
longer chain which will contradict with the maximality of $M$. So,
$y$ must lie already in $M$. Thus, $M \in \M _y $.
\end{proof}

We have the following:

\begin{prop}\label{prop:posetprimreduction}
Let $P$ be a poset and let $x,y \in P $ such that $x \prec y$. Then,
$\D(P) \searrow _{\LC} \D (P \backslash \{ x \}) $.
\end{prop}

\begin{proof} Consider the primitive linear coloring $\ka$ that involves
only $x$ and $y$. The proposition follows from the fact that $\D
(P)_k =\del _{\D (P)} (x)=\D (P \backslash \{ x \})$.
\end{proof}

It is easy to see that if an element is minimal or maximal, then it
dominates all other elements. So, if a poset has a minimal or
maximal element, then it is $\LC$-contractible.

Now, we consider monotone poset maps and prove a reduction theorem
for them.

\begin{defn}\label{defn:monotonemap}
Let $P$ be a poset. An order-preserving map $\varphi \colon P\to P$
is called a \emph{monotone map} if either $x\leq \varphi(x)$ or
$x\geq \varphi(x)$ for any $x\in P$. If $\varphi$ is a monotone map
which also satisfies $\varphi ^2=\varphi$, then it is called a
\emph{closure operator} on $P$.
\end{defn}

Note that when $\varphi : P \to P$ is a closure operator then $\Fix
(\varphi )=\varphi (P)$, and the equality $P=\varphi (P)$ holds only
when $\varphi$ is the identity map.

\begin{lem}\label{lem:monotonemap} Let $P$ be a finite poset, and let
$\psi : P \to P$ be a monotone map on $P$ which is different than
the identity map. Then there exists a $x \in P \backslash \Fix
(\psi)$ such that $x \prec \psi (x) $.
\end{lem}

\begin{proof}  Assume to the contrary that for all $x \in P \backslash
\Fix(\psi)$, we have $x \not \prec \psi (x)$. Start with $y_0
\in P\backslash \Fix(\psi)$ such that $y_0 \not \prec \psi
(y_0 )$. This means that there exists an element $y_1 \in P$ such
that $y_1$ is comparable with $y_0$ but not with $\psi (y_0)$.

Note that since $\psi $ is a monotone map either $y_0 < \psi
(y_0)$ or $\psi (y_0)< y_0$ holds. We look at each case
separately.

Case 1: Assume $y_0 < \psi (y_0)$ holds. Then, we must have $y_0
<y_1$, because otherwise we have $y_1 < y_0 < \psi (y_0)$ which
contradicts the assumption that $y_1$ and $\psi (y_0)$ are not
comparable. Also note that $y_1 $ cannot be an element of $\Fix
(\psi )$, because otherwise $y_1=\psi (y_1) < \psi (y_0)$ implies
that $y_1$ and $\psi (y_0)$ are comparable, which is again a
contradiction. So, we have $y_1 \in P \backslash \Fix(\psi)$.

Now, let's apply the same arguments for $y_1$. First we have $y_1
\not \prec \psi (y_1)$ by our starting assumption, so there exists a
$y_2$ such that $y_2$ comparable with $y_1$ but not with $\psi
(y_1)$. Since $\psi $ is a monotone map, we again have either $y_1 <
\psi (y_1)$ or $\psi (y_1)< y_1$. Now we claim that actually the
second inequality cannot hold. Suppose it holds, i.e., $\psi (y_1)<
y_1$. Then we get $\psi (y_0)< \psi (y_1 ) <y_1$ which gives
$\psi(y_0)$ and $y_1$ are comparable and hence a contradiction. So,
we have $y_1 < \psi (y_1)$. This allows us to continue in the same
way and obtain an infinite ascending sequence $y_0 < y_1 < y_2 <
\cdots $ of distinct elements in $P$. But, this is in contradiction
with the fact that $P$ is a finite poset.

Case 2:  Assume $y_0 > \psi ( y_0)$ holds. Then, arguing as above we
find a descending infinite sequence $y_0 > y_1 > y_2 > \cdots $ of
distinct elements in $P$ and again reach a contradiction.
\end{proof}

The main result of this section is the following:

\begin{thm}\label{thm:monotonemap}
Let $\varphi\colon P\to P$ be a closure operator on a finite poset
$P$. Then, $\D(P)\searrow_{\LC} \D(\varphi (P) )$.
\end{thm}

\begin{proof}
We will prove the result by induction on $n=|P \backslash \varphi
(P)|$. If $n=0$, then there is nothing to prove. So assume $n \geq
1$, i.e., $\varphi$ is not identity. Then, by Lemma
\ref{lem:monotonemap} there exists a $x \in P\backslash \varphi (P)$
such that $x \prec \varphi (x)$. By Proposition
\ref{prop:posetprimreduction}, we have $\D (P ) \searrow _{\LC} \D
(P \backslash \{ x \} )$. Since $x \not \in \varphi (P)$, the
restriction of $\varphi$ to $P \backslash \{ x \}$ induces a closure
operator $\overline{\varphi} : P \backslash \{ x\} \to P \backslash
\{ x \}$. Applying the induction assumption, we obtain $\D (P
\backslash \{ x\} ) \searrow_{\LC} \D(\overline{\varphi}
(P\backslash \{ x\}))$ which gives $\D (P \backslash \{ x\} )
\searrow_{\LC} \D(\varphi (P))$ since $\overline{\varphi} (P
\backslash \{x\})=\varphi(P)$. Combining this with the above
reduction, we conclude that $\D (P) \searrow _{\LC} \D( \varphi
(P))$.
\end{proof}

\begin{cor}\label{cor:coatom}
For a finite poset $P$, if $\bar{x}=\bigwedge \{c\in co(P)\colon
x\leq c\}$ exists for all $x\in P$ then $P\searrow_{\LC} R$, where
$R=\{\bar{x}\ |\ x\in P\}$. If in addition, $\bigwedge co(P)$ exists
then  $\D (P) $ is $\LC$-contractible.
\end{cor}

\begin{proof}
The map $\varphi \colon P\to P$ defined by $\varphi(x)=\bar{x}$ is a
closure operator. Hence, by Theorem \ref{thm:monotonemap},
$\D(P)\searrow_{\LC} \D(R)$, since $\Fix (\varphi )=\varphi(P)=R$.
On the other hand, when it exists, $\bigwedge co(L)$ is the minimal
element of $R$, therefore $\D(R)$ is $\LC$-contractible so is
$\D(P)$.
\end{proof}

In particular, the above corollary says that the proper part of a
lattice is $\LC$-reducible to the proper part of the sublattice of
elements that are the meet of coatoms. This result is well-known
when the $\LC$-reduction is replaced by homotopy equivalence (see
Theorem 10.8 in \cite{Bj}).

Another interesting invariant in poset theory is the order dimension
of a poset which is defined as follows:

\begin{defn}\label{defn:orddim} The \emph{order dimension} of a
finite poset $P$, denoted by $\ordim (P)$, is defined to be the
smallest integer $n$ such that $P$ can be embedded in $\N^n$ as an
induced subposet (an induced subposet is a subposet which inherits
all the relations of the poset.)
\end{defn}

There is a very nice paper by Reiner and Welker \cite{RW} which
proves that the order dimension of a lattice $L$ is greater that
$\homdim ( \overline{L} ) +2$ where $\overline{L}$ denotes the
proper part of the lattice $L$. Recall that there is a similar
inequality for the linear dimension of a poset (see Proposition
\ref{prop:lhomdim}). The obvious question is whether there is any
connection between the order dimension of a lattice and the linear
dimension of its proper part. Unfortunately these invariants are not
comparable by inclusion as the following examples show.

\begin{exmp}\label{ex:ordimversuslindim}
Consider the poset $P$ which is an antichain with three elements.
Let $L$ be the lattice obtained form $P$ by adding minimal and
maximal elements. It is clear that $\overline{L}=P$ has linear
dimension exactly $3$. But, the order dimension of $L$ is equal to
$2$ since we can embed $L$ in $\N^2$ by taking the minimal element
to $(0,0)$, the maximal element to $(2, 2)$ and the $3$ middle
points to the points $(0, 2), (1,1), (2,0)$. This shows that there
is a lattice $L$ where $\ordim (L)< \lindim (\overline{L})$.

For the other direction, consider the poset $P=\{ a,b,c \}$ where $a
\leq b$, $a \leq c$, and $b$ and $c$ are not comparable. It is easy
to see that $P$ is $\LC$-reducible to a point so $\lindim(P)=1$. Let
$L$ be the lattice obtained from $P$ by adding $\hat 0$ and $\hat
1$. It is clear that $L$ is not linear, so $\ordim (L)>1=\lindim
(\overline{L})$.
\end{exmp}

We end the section with an application of  Corollary
\ref{cor:coatom} to subgroup lattices.

\begin{cor}\label{cor:p-group}
Let $G$ be a finite $p$-group ($p$ a prime). Then,
$\overline{\Lattice(G)}$  is $\LC$-contractible if and only if $G$
is not elementary abelian, where $\Lattice(G)$ is the subgroup
lattice of $G$.
\end{cor}

\begin{proof}
It is known that if $G$ is elementary abelian, then the Euler
characteristic of $\overline{\Lattice(G)}$ is bigger than $1$. Thus,
$\overline{\Lattice(G)}$ cannot be $\LC$-contractible. Conversely,
if $G$ is not elementary abelian, then the intersection of the
maximal subgroups of $G$ is non-trivial. Therefore, by Lemma
\ref{cor:coatom}, $\overline{\Lattice(G)}$ is $\LC$-contractible.
\end{proof}

%%%%%%%%%%%%%%%%%%%%%%%%%%%%%%%%%%%%%%%%%%%%%%%%%%%%%%%%%%%%%%%%%%%%

\section{Linear graph colorings}\label{sect:graphcoloring}

In this final section, we consider linear colorings of neighborhood
complexes associated to simple graphs.

Let $G=(V,E)$ be a simple graph. We recall that a (vertex) coloring
of $G$ is a surjective mapping $\nu\colon V\to [n]$ such that
$\nu(x)\neq \nu(y)$ whenever $(x,y)\in E$. The \emph{neighborhood
complex} of $G$, denoted by $\Neigh(G)$, is defined as the
simplicial complex whose simplices are those subsets of $V$ which
have a common neighbor. We start with the following easy
observation.

\begin{prop}
\label{prop:coloringN(G)} Let $G=(V,E)$ be a simple graph and let
$\Neigh (G)$ denote its neighborhood complex. If $\ka \colon V\to
[k]$ is a $k$-linear coloring of $\Neigh(G)$, then $\ka$ is a
coloring of the underlying graph $G$.
\end{prop}

\begin{proof} Assume that $\ka$ is not a coloring of the underlying
graph $G$. Therefore, there exist $x, y \in V$ such that $(x, y)\in
E$ and $\ka(x)=\ka(y)$. By the definition of a linear coloring,
either $\F (x) \subseteq \F (y)$ or $\F(y) \subseteq \F (x)$. So,
without loss of generality, assume $\F(x) \subseteq \F (y)$. Let
$\Neigh (z)$ be a facet of $\Neigh (G)$ such that  $\Neigh (y)
\subseteq \Neigh (z)$. Since there is an edge between $x$ and $y$,
we have $x\in \Neigh (y)$, and hence $x \in \Neigh (z)$. This
implies that $\Neigh (z)\in \F (x)$, and gives $\Neigh (z) \in \F
(y)$. Therefore, $y \in \Neigh (z)$ and hence $z \in \Neigh (y)$.
However, together with $\Neigh (y) \subseteq \Neigh (z)$, this
implies $z \in \Neigh (z)$ which is a contradiction since $G$ is a
simple graph and has no loops.
\end{proof}

The following is immediate:

\begin{cor}\label{cor:graphchromatic}
For any graph $G$, we have $\lchr(\Neigh(G))\geq \chi(G)$, where
$\chi(G)$ denotes the (vertex) chromatic number of $G$.
\end{cor}

It is easy to see that a coloring of $G$ may not give rise to a
linear coloring of its neighborhood complex $\Neigh(G)$. So, in
general the equality does not hold.

\begin{exmp}\label{ex:hexagon}
Consider the graph which is an hexagon, i.e., $G=(V, E)$ with
$V=\{v_1, \dots, v_6\}$ and $E=\{ (v_i, v_{i+1} ) \ | \ 1\leq i \leq
5 \} \cup \{ (6,1) \} $. Note that $\chi (G)=2$, but $
\lchr(\Neigh(G))=6$ since $\Neigh (G)$ is a disjoint union of two
triangles.
\end{exmp}

We now give a sufficient condition for a coloring of a graph to be a
linear coloring of its neighborhood complex.

\begin{prop}\label{prop:Ncoloring}
A coloring $\nu\colon V\to [k]$ of $G=(V,E)$ is a $k$-linear
coloring of $\Neigh(G)$ if either $\Neigh(v)\subseteq \Neigh(u)$ or
$\Neigh(u)\subseteq \Neigh(v)$ holds for every $x,y\in V$ with
$\nu(x)=\nu(y)$.
\end{prop}

\begin{figure}[ht]
\begin{center}
\psfrag{1}{$1$}
\psfrag{2}{$2$}\psfrag{3}{$3$}\psfrag{u}{$u$}\psfrag{v}{$v$}
\psfrag{w}{$w$}\psfrag{x}{$x$}\psfrag{y}{$y$}\psfrag{z}{$z$}
\psfrag{G}{$G$}\psfrag{N}{$\Neigh(G)$}
\includegraphics[width=4in,height=1.2in]{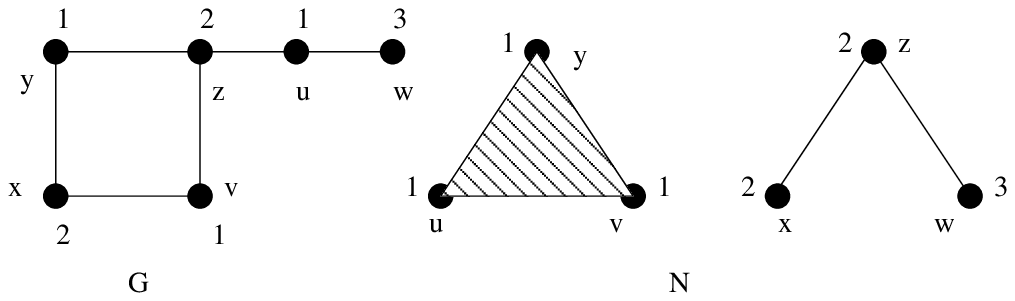}
\end{center}
\caption{}\label{colG2}
\end{figure}

\begin{proof}
Assume that whenever $\nu(u)=\nu(v)$ for any two vertices $u, v\in
V(G)$, then one of the inclusions $\Neigh(v)\subseteq \Neigh(u)$ or
$\Neigh(u)\subseteq \Neigh(v)$ holds. Let $u, v\in V(G)$ be two such
vertices and let $\Neigh(u)\subseteq \Neigh(v)$. To verify that
$\F(u)\subseteq \F(v)$, let $\Neigh(y)$ be a facet of $\Neigh(G)$
containing $u$. Then we must have $y\in \Neigh(v)$, since $y\in
\Neigh(u)\subseteq \Neigh(v)$. Hence, $v\in \Neigh(y)$.
\end{proof}
The converse of Proposition \ref{prop:Ncoloring} does not hold in
general as illustrated in Figure \ref{colG2}. It is easy to see that
the given vertex coloring of $G$ is indeed a linear coloring of
$\Neigh(G)$ with $\nu(u)=\nu(v)=1$; however, there is no inclusion
relation between the neighborhoods of $u$ and $v$.


\begin{thebibliography}{00}

\bibitem{Bj} A.~Bj\"{o}rner,
\newblock{\em Topological Methods,}
\newblock{Handbook of Combinatorics, R.~Graham, M.~Gr\" otscel,
and~L. Lov\' asz (eds),}
\newblock{North-Holland/Elsevier, Amsterdam, 1995, 1819-1872.}

\bibitem{BV} A.~Bj\"{o}rner and S.~Vre\` cica,
\newblock{\em On $f$-vectors and Betti numbers of multicomplexes,}
\newblock{Combinatorica, {\bf 17} (1997), 53-65.}

\bibitem{KSS} J.~Kahn, M.~Saks, and D.~Sturtevant,
\newblock{\em A topological approach to evasiveness,}
\newblock{Combinatorica, {\bf 4} (1984), 297--306.}

\bibitem{Ko1} D.N.~Kozlov,
\newblock{\em Collapsing along monotone poset maps,}
\newblock{to appear in the International Journal of Mathematics and
Mathematical Sciences,}
\newblock{arXiv:math.CO/0503416.}

\bibitem{Ko2} D.N.~Kozlov,
\newblock{\em Simple homotopy types of Hom-complexes, neighborhood complexes,
Lov\'asz complexes, and atom crosscut complexes,}
\newblock{to appear in Topology Appl.,}
\newblock{arXiv:math.AT/0503613.}

\bibitem{Ku} H.~Kurzweil,
\newblock{\em A combinatorial technique for simplicial complexes
and some applications to finite groups,}
\newblock{Discrete Math., {\bf 82} (1990), 263-278.}

\bibitem{Mu} J.R.~Munkres,
\newblock{\em Elements of Algebraic Topology,}
\newblock{Addison-Wesley Pub., New York, 1993.}

\bibitem{Lo} L.~Lov\'{a}sz,
\newblock{\em Kneser's conjecture, chromatic number, and homotopy,}
\newblock{J. Combinatorial Theory, Series {\bf A}, {\bf 25} (1978), 319-324.}

\bibitem{Qu} D.~Quillen,
\newblock{\em Homotopy properties of the poset of nontrivial
$p$-subgroups of a group,}
\newblock{Adv. in Math., {\bf 28(2)} (1978), 101-128.}

\bibitem{RW} V.~Reiner and V.~Welker,
\newblock{\em A homological lower bound for order dimension of lattices,}
\newblock{Order, {\bf 16} (1999),165-170.}

\bibitem{S1} R.P.~Stanley,
\newblock{\em Combinatorics and Commutative Algebra,}
\newblock{Progress in Mathematics, {\bf 41}, Birkh\"{a}user, Boston, 1997.}

\bibitem{We} V.~Welker,
\newblock{\em Constructions preserving evasivenes and collapsibility,}
\newblock{Discrete Math., {\bf 207} (1999), 243-255.}


\end{thebibliography}
\end{document}